\documentclass{amsart}

\usepackage{amsmath,amssymb,amsthm}
\usepackage{graphicx,tikz}
\newtheorem{theorem}{Theorem}
\newtheorem*{thm}{Theorem}

\theoremstyle{definition}

\theoremstyle{remark}

\begin{document}

\title[]{On Concavity of solutions of \\the nonlinear Poisson Equation}
\thanks{S.S. is supported by the NSF (DMS-2123224) and the Alfred P. Sloan Foundation.}

\author[]{Stefan Steinerberger}
\address{Department of Mathematics, University of Washington, Seattle, WA 98195, USA}
\email{steinerb@uw.edu}

\begin{abstract} We consider the nonlinear Poisson equation $-\Delta u = f(u)$
 in domains $\Omega \subset \mathbb{R}^n$ with Dirichlet boundary conditions on $\partial \Omega$. We show (for monotonically increasing concave $f$ with small Lipschitz constant) that if $D^2 u$ is negative semi-definite on the boundary, then $u$ is concave. A conjecture of Saint Venant from 1856 (proven by Polya in 1948) is that among all domains $\Omega$ of fixed measure, the solution of $-\Delta u =1$ assumes its largest maximum when $\Omega$ is a ball. We extend this to $-\Delta u =f(u)$ for monotonically increasing $f$ with small Lipschitz constant. \end{abstract}

\maketitle

\section{Introduction}
\subsection{Introduction}
We study the behavior of solutions of 
\begin{align*}
-\Delta u &= f(u) \qquad \mbox{in}~\Omega\\
u &= 0 \qquad \quad \mbox{on}~\partial \Omega,
\end{align*}
where $\Omega \subset \mathbb{R}^n$ is assumed to be a set with smooth boundary and $f:\mathbb{R}_{\geq 0} \rightarrow \mathbb{R}_{\geq 0}$ is assumed to be a nonnegative function with $f(0) > 0$. One natural assumption is that solutions of this equation should `inherit' the simplicity of the underlying domain and also be `simple'.  An optimistic guess is that such functions are concave but this is already violated for very simple equations such as $-\Delta u = 1$ (see Fig. 1). Another simple property that one might hope for is perhaps that the level sets are convex (conjectured by P.-L. Lions \cite{lions} in 1981): this was disproven only rather recently
by Hamel, Nadirashvili \& Sire \cite{hamel} in 2016.

\begin{center}
\begin{figure}[h!]
\includegraphics[width =0.45\textwidth]{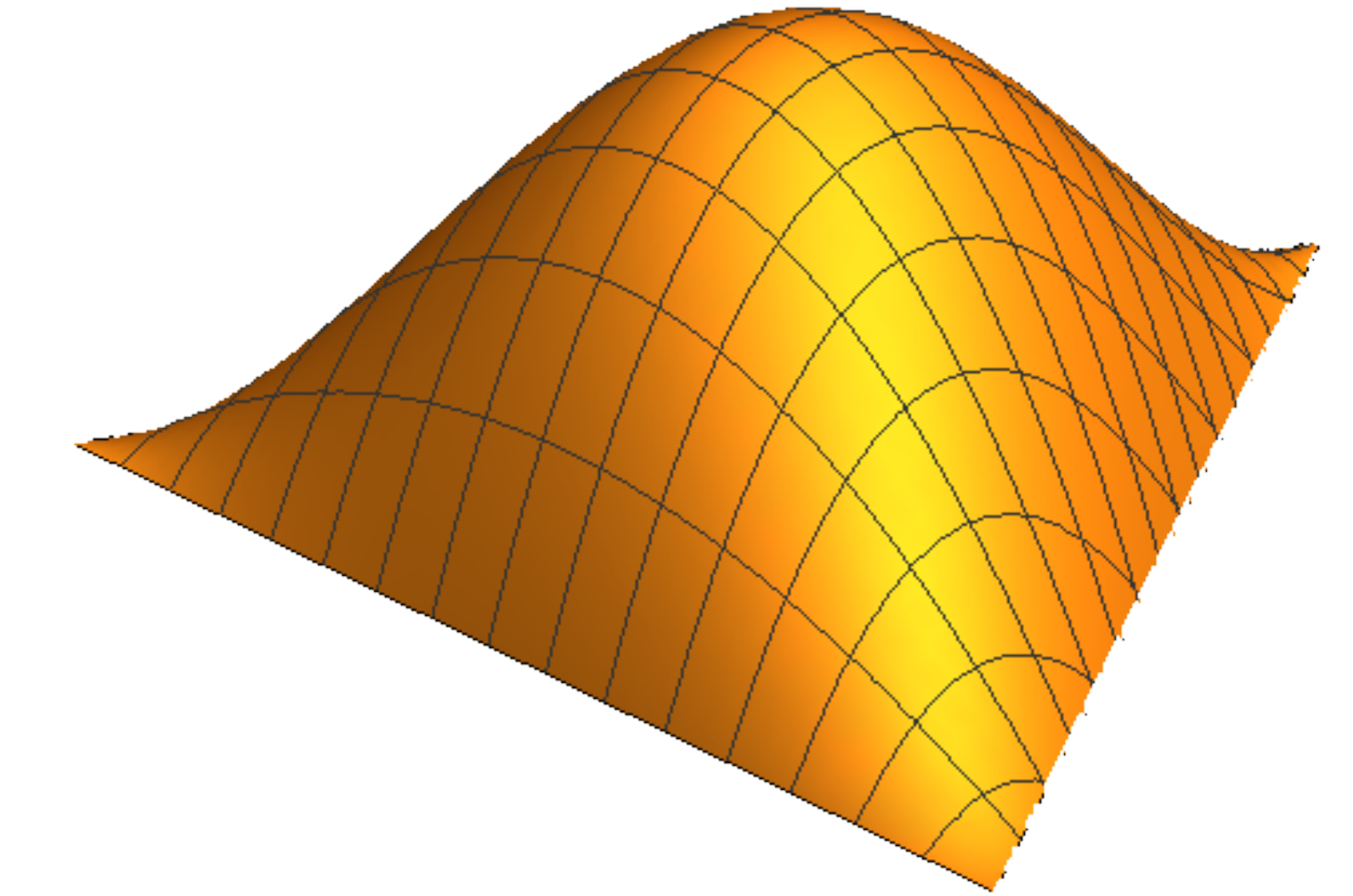}
\caption{$-\Delta u = 1$ on a triangle need not be concave (see \cite{keady}).}
\end{figure}
\end{center} 

 Nonetheless, it is known that the level sets are convex for some special cases: Makar-Limanov \cite{makar} proved that solutions of $-\Delta u = 1$ have the property that $\sqrt{u}$ is convex (and thus $\left\{x \in \Omega: u(x) \geq t\right\}$ is a convex set). Brascamp \& Lieb \cite{brascamp} showed that if $u$ is the first Laplacian eigenfunction on $\Omega$, i.e. the solution of $-\Delta u = \lambda_1(\Omega) u$, then $\log{u}$ is concave. Keady \cite{keady0} proved that for
$$f(u) = u^{\gamma} \quad \mbox{and} \quad 0 < \gamma < 1,$$
 the function $u^{(1-\gamma)/2}$ is concave, this was extended by Kennington \cite{ken} to higher dimensions. There are many results of such a flavor: we refer to 
\cite{beck2, henrot, hoskins, keady, lu, stv1, makar, stein} for solutions of the torsion problem $-\Delta u = 1$ alone and to 
\cite{arango, bandle, beck1, beck3, bian, borell, borell2,  caf, caf2, cha, chau1, finn, gab, gugg, hen, kawohl0, kor, kor2, ma} for general equations.

\subsection{The Problem.}
We were motivated by two simple arguments. First, if $\Omega$ is a disk, then solutions of such equations tend to be radially symmetric (this is very well understood \cite{gidas}). Secondly, there is a simple reason why convexity of $\Omega$ might not be `right' condition. Consider the equation
$ -\Delta u = 1$ which can be non-concave close to corners (see Fig. 1).
One would expect something similar to be true for $-\Delta u = f(u)$ in general:  when $x$ is close to the boundary, then $u$ is close to 0 because of the Dirichlet boundary conditions and $f(u(x)) \sim f(0)$ and the solution should behave, at least locally, a bit like a multiple of $-\Delta u = 1$. Convexity of the domain $\Omega$ is a natural assumption but it should be equally important to rule out corners and, by approximation arguments, regions of the boundary with large curvature.
 This motivates the following (intentionally vague) 
\begin{quote} 
\textbf{Question.} Under what conditions on $f>0$ (`sufficiently nice') and under what conditions on the convex domain $\Omega \subset \mathbb{R}^n$ (sufficiently `round') is it true that solutions of
\begin{align*}
-\Delta u &= f(u) \qquad \mbox{in}~\Omega\\
u &= 0 \qquad \quad \mbox{on}~\partial \Omega,
\end{align*}
are concave everywhere in $\Omega$?
\end{quote}

We are not aware of any results of this flavor with the exception of one:  Kosmodemyanskii \cite{kos0} showed that if every parabola osculating with the boundary of $\Omega \subset \mathbb{R}^2$ contains $\Omega$ in its interior, then the solution of $-\Delta u = 1$ is concave.

\section{The Results} 
\subsection{Propagating Concavity.}
 We start, for illustrative purposes, by considering a simple equation which has been studied for a long time
\begin{align*}
-\Delta u &= 1 \qquad \mbox{in}~\Omega \\
u &= 0 \qquad  \mbox{on}~\partial \Omega.
\end{align*}
It is known that if $D^2 u$ is negative semi-definite on the boundary, then $u$ is concave inside $\Omega$ (a `propagation of concavity' phenomenon). We will illustrate our main idea by giving a new proof of this result (our general approach becomes particularly simple when $f \equiv 1$ which is thus well suited for outlining the idea).

\begin{thm}[Keady \& McNabb, \cite{keady}]
Let $\Omega \subset \mathbb{R}^2$ and let $u$ be the solution of $-\Delta u = 1$ with Dirichlet boundary conditions. If $D^2u$ is negative semi-definite on $\partial \Omega$,
then $u$ is concave on all of $\Omega$.
\end{thm}
We emphasize that our approach is quite different and uses a stochastic representation of the solution together with a bootstrapping argument.
Our first main result is a generalization of this result to equations of the type $-\Delta u =f(u)$.

\begin{theorem}[Propagating Concavity from the Boundary]
Let $\Omega \subset \mathbb{R}^n$ have a smooth boundary, let $-\Delta u = f(u)$ and assume $f > 0$, $f' \geq 0$ and $f'' \leq 0$ and
 $$  f'(0) \leq \frac{2 n \omega_n^{2/n}}{|\Omega|^{2/n}},$$
 where $\omega_n$ is the volume of the unit ball in $\mathbb{R}^n$.
If $D^2 u$ is negative semi-definite for all $x \in \partial \Omega$, then $u$ is concave in $\Omega$. 
\end{theorem}
In the special case of the torsion function $f \equiv 1$, the condition is trivially satisfied independently of the domain. One of the reasons the concavity of such solutions is difficult to quantify is that sometimes they are `barely' concave. Consider the solution of $-\Delta u = 1$ in a long rectangle. An explicit computation (see \cite{stein}) shows that if $x_0$ is the point of the rectangle in which $u$ assumes its maximum, then the largest eigenvalue of the Hessian satisfies, for some $c_1, c_2 > 0$,
$$ \lambda_{\max}(D^2 u(x_0)) \sim - c_1 \exp\left( - c_2 \frac{\mbox{diam}(\Omega)}{\mbox{inrad}(\Omega)} \right).$$
As the rectangle becomes more eccentric, this eigenvalue is exponentially close to 0: the function is concave but it becomes
harder to tell from the behavior of the eigenvalues of the Hessian. Our argument will recover this type of
phenomenon and lead to a new explanation why this is the case (see \S 3.2).

\subsection{A Rearrangement Result.}
Upper bounds on $\| u\|_{L^{\infty}(\Omega)}$ are a classical subject. It is generally understood that symmetric
decreasing rearrangement is a process that tends to increase $L^p$ norms of solutions.
The most fundamental formulation of this statement is perhaps Talenti's rearrangement principle \cite{talenti, talenti2}.
If
$$  - \Delta u = f \geq 0 \qquad \mbox{in}~\Omega$$
with Dirichlet boundary conditions on $\partial \Omega$, then solving the same problem with $f$ rearranged on a ball $B$
with the same measure as $\Omega$, i.e.
$$ -\Delta v = f^* \qquad \mbox{in}~B,$$
leads to a solution for which $v \geq u^*$ in a pointwise sense. Since $u^*(0) = \|u\|_{L^{\infty}(\Omega)}$, we
conclude that $v$ has a larger maximum.
One could wonder whether there is a similar result for $-\Delta u = f(u)$ but little seems to be known.
 The book of Baernstein \cite{baern} discusses the case of
decreasing $f$ following a result of Weitsman \cite{weitsman}.
\begin{theorem}[Rearrangement Theorem] Let $\Omega \subset \mathbb{R}^n$ be a domain with smooth boundary, let $f > 0$, $f' \geq 0$ and
 $$ \max_{t> 0} f'(t)  <  \frac{2n \omega_n^{2/n}}{ |\Omega|^{\frac{2}{n}}},$$
 where $\omega_n$ is the volume of the unit ball in $\mathbb{R}^n$.
The solution of
\begin{align*}
-\Delta u &= f(u) \qquad \mbox{in}~\Omega\\
u &= 0 \qquad \quad \mbox{on}~\partial \Omega,
\end{align*}
satisfies $u^* \leq v$, where $v$ is the solution on the ball with same measure as $\Omega$. 
\end{theorem}
We note that in the case of $f \equiv 1$, this result was first suggested by Saint Venant in 1856. The first proof was given by Polya \cite{pol} in 1948.  
Nowadays, Talenti's rearrangement principle \cite{talenti, talenti2} can be used to give a very simple proof. Theorem 2 can be considered a nonlinear analogue of these ideas.  It is not clear to us whether the upper bound on $f'$ is necessary.

\section{Proof of Theorem 1}
This section gives a proof of Theorem 1.
We first discuss the case of $-\Delta u = 1$ which is particularly well suited to illustrate the main ideas behind the argument. 
\subsection{An Identity.} We start with the following identity. If
$$ -\Delta u = 1 \qquad \mbox{in}~\Omega$$
and $u=0$ on the boundary, then for any $x \in \Omega$ and any unit vector $\|n\|=1$
$$ \frac{\partial^2 u}{\partial n^2}(x) = \int_{\partial \Omega}  \frac{\partial^2 u}{\partial n^2}(y)~ d\omega_x(y),$$
where $\omega_x$ is the harmonic measure on the boundary $\partial \Omega$ induced by the point $x \in \Omega$.
In particular, it suffices to prove that $D^2 u$ is negative definite on the boundary to deduce global concavity. Theorem 2 will
establish a more general form of the identity. Since the case $f \equiv 1$ is particularly simple,
we give an independent and much simpler argument here. In particular, the proof uses so little that it also naturally extends to higher derivatives.

\begin{center}
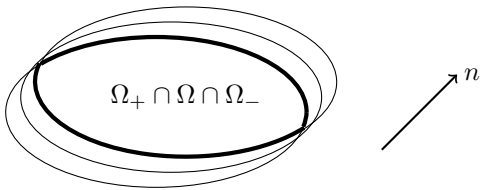
\begin{figure}[h!]
\begin{tikzpicture}
 \draw[] (0,0) arc (0:360:2cm and 1cm);
   \draw[ultra thick] (0,0) arc (359:348:2cm and 1cm);
  \draw[ultra thick] (0,0) arc (0:140:2cm and 1cm);
  \draw[] (0.2,0.2) arc (0:360:2cm and 1cm);
 \draw[] (0.4,0.4) arc (0:360:2cm and 1cm);
  \draw[ultra thick] (0.4-3.95,0.65) arc (166:322:2cm and 1cm);
  \draw[thick, ->] (1,-0.5) -- (2,0.5);
  \node at (2.2, 0.5) {$n$};
  \node at (-1.6, 0.2) {$ \Omega_{+} \cap \Omega \cap \Omega_{-}$};
\end{tikzpicture}
\caption{Shifting the domain by $\pm \varepsilon$ in direction $n$.}
\end{figure}
\end{center}
\begin{proof}[Proof of the Identity.] The idea is to use the differential quotient
$$ \frac{\partial^2 u}{\partial n^2}(x) = \lim_{\varepsilon \rightarrow 0} \frac{u(x+\varepsilon n) - 2 u(x) + u(x- \varepsilon n)}{\varepsilon^2}.$$
Instead of interpreting this as one function evaluated at three nearby points, we will consider it as the evaluation of three different functions in the
same point. The two additional functions are the solutions of the same PDE in the shifted domains
$$ -\Delta u = 1 \qquad \mbox{in}~\Omega \pm \varepsilon n.$$

We denote the three functions by $u_{+}, u, u_{-}$. Naturally, since
$$ \Delta \left( \frac{u_{+} - 2 u + u_{-}}{\varepsilon^2} \right) = 0 \qquad \mbox{in}~ \Omega_{+} \cap \Omega \cap \Omega_{-},$$
the differential quotient is a harmonic function in the intersection of the three domains. We know how to solve Laplace's equation
inside a domain, the solution is given by integrating harmonic measure against boundary data. It thus remains to determine
boundary data but this is simply given by the second directional derivative at the boundary in direction $n$. Letting $\varepsilon \rightarrow 0$
leads to the result.
\end{proof}

This identity immediately implies Theorem 1 for the special case $f \equiv 1$.

\subsection{A Remark on Eccentric Level Sets.} The author proved \cite{stein} that if $\Omega \subset \mathbb{R}^2$ is a convex domain, then the (unique) maximum of the solution of
$$ -\Delta u = 1 \qquad \mbox{in}~\Omega$$
has the property that the determinant of its Hessian cannot be too small: level sets close to the maximum of the solution are
asymptotically ellipses whose eccentricity can be bounded in terms of the geometry of the domain: for some universal $c_1, c_2 > 0$, the largest
eigenvalue of the Hessian on the maximum $x_0 \in \Omega$ satisfies
$$ \lambda_{\max}(D^2 u (x_0)) \leq - c_1 \exp\left( - c_2 \frac{\mbox{diam}(\Omega)}{\mbox{inrad}(\Omega)}\right).$$

The proof uses conformal mapping and only works in two dimensions. The formula
$$ \frac{\partial^2 u}{\partial n^2}(x) = \int_{\partial \Omega}  \frac{\partial^2}{\partial u^2}(y)~ d\omega_x(y)$$
is more generally applicable and
may be helpful in explaining why it is possible that the level sets close to the maximum are so remarkably eccentric. 
In long domains,
the harmonic measure is rather small (exponentially in eccentricity) at the endpoints. 
\begin{center}
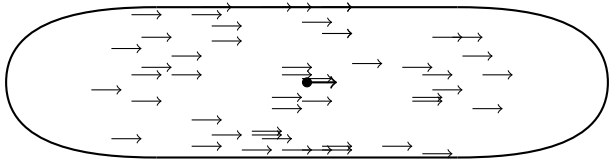
\begin{figure}[h!]
\begin{tikzpicture}[scale=2]
\draw [thick] (0,0) -- (2,0);
\draw [thick] (2,0) to[out=0, in =270] (3, 0.5) to[out=90, in =0] (2,1);
\draw[thick] (2,1) -- (0,1);
\draw [thick] (0,1) to[out=180, in=90] (-1, 0.5) to[out=270, in=180] (0,0);
\filldraw (1, 0.5) circle (0.03cm);
\draw [thick, ->] (1, 0.5) -- (1.2, 0.5);

   \foreach \y in {0,...,50}
      {
 \pgfmathsetmacro{\Xa}{random(40)/15-0.5}
  \pgfmathsetmacro{\Xb}{random(40)/40}
\draw [->] (\Xa, \Xb) -- (\Xa + 0.2, \Xb);
      }
 \end{tikzpicture}
\caption{An example domain: most of the randomly traveling vectors hit the boundary along the parallel segment.}
\end{figure}
\end{center}
However, if $n$ is chosen to be parallel to the long
direction, then most contributions in the integral are going to be very small (see Fig. 4). 
It is conceivable that the representation formula in \S 3.3 could be used to prove that the same phenomenon
also occurs for $- \Delta u = f(u)$ for suitable $f$ and to provide bounds for the largest eigenvalue of the Hessian; this
might be an interesting avenue for future research.

\subsection{A Representation Formula} 
Let $\Omega \subset \mathbb{R}^n$ and let $u: \Omega \rightarrow \mathbb{R}$ solve
\begin{align*}
-\Delta u &= f(u) \qquad \mbox{in}~\Omega\\
u &= 0 \qquad \quad \mbox{on}~\partial \Omega,
\end{align*}
for any $x \in \Omega$ and any unit length vector $\|n\|=1$, the identity
\begin{align*}
\frac{\partial^2 u}{\partial n^2}(x) &= \mathbb{E}~ \frac{1}{2}   \int_0^{\infty}   f''(u(\omega_x(s))) \left( \frac{\partial u}{\partial n}(\omega_x(s)) \right)^2 +  f'(u(\omega_x(s))) \frac{\partial^2 u}{\partial n^2}(\omega_x(s)) ds \\
&+ \int_{\partial \Omega} \frac{\partial^2 u}{\partial n^2}(y) d\omega_x(y),
\end{align*}
where $\omega_x(s)$ denotes a Brownian motion started in $x$ after $s$ units time with the convention that it gets absorbed by the boundary $\partial \Omega$ and stops there (in particular, the first integral could have also been written with the upper limit $\infty \wedge \tau$, where $\tau$ is the first hitting time). In the second integral, $d \omega_x$ denotes the harmonic measure induced by $x$ on the boundary $\partial \Omega$. We use a stochastic representation of the solution -- since this may be hard to find in the literature and since precise constants matter for our approach, we quickly discuss a self-contained heuristic derivation. Let $u:\mathbb{R}^n \rightarrow \mathbb{R}$ be a smooth function and suppose $\omega(t)$ is an $n$-dimensional Brownian motion started in the origin. Each coordinate of $\omega(t)$ is a one-dimensional Brownian motion. In particular, we have
$$ \mathbb{E}~ \| \omega(t)\|^2 = \mathbb{E} ~\omega_1(t)^2 + \dots + \mathbb{E}~\omega_n(t)^2 = n t.$$
We will now compute $\mathbb{E}~u(\omega(t))$ for small values of $t$. Using a Taylor expansion of $u$, we have up to lower order terms
$$ u(x) = u(0) + \left\langle \nabla u, x\right\rangle +   \frac{1}{2} \left\langle x, (D^2 u(0)) x\right\rangle + \mbox{l.o.t.}$$
It remains to set $x = \omega(t)$ and take an expectation with respect to $\omega$.
The symmetry of Brownian motion guarantees that the linear/gradient term does not play any role. It remains to understand the quadratic term
$$ \mathbb{E} ~u(\omega(t)) = u(0) + \frac{1}{2} \cdot \mathbb{E}  \left\langle \omega(t), (D^2 u(0)) \omega(t)\right\rangle + \mbox{l.o.t.}$$
Brownian motion is radially distributed around its point of origin. In terms of probability, we can thus write $\omega(t) = \| \omega(t)\| \cdot X,$ where $X \sim U(\mathbb{S}^{n-1})$ is a uniformly chosen random point on the unit sphere. Then, since radial and angular distribution of a Gaussian can be decoupled,
\begin{align*}
 \mathbb{E}_{\omega}  \left\langle \omega(t), (D^2 u(0)) \omega(t)\right\rangle &= \mathbb{E}_{\omega} \left( \| \omega(t)\|^2 \cdot    \left\langle X, (D^2 u(0)) X\right\rangle \right)\\
 &= \left( \mathbb{E}_{\omega} \| \omega(t)\|^2 \right)  \cdot   \left( \mathbb{E}_{X} \left\langle X, (D^2 u(0)) X\right\rangle \right).
\end{align*}
The first of these terms was already computed above. As for the second quantity, using symmetry of $X$to expand $D^2(u(0))$ into its eigenvectors and eigenvalues. Then
$$ \mathbb{E}_{X} \left\langle X, (D^2 u(0)) X\right\rangle = \mathbb{E}_X \sum_{k=1}^{n} \lambda_k \left\langle v_k, X\right\rangle^2.$$
However, the distribution of $X$ is rotationally invariant and thus $\mathbb{E}_X \left\langle v, X \right\rangle^2$ is independent of the vector $v$ (as long as $v$ has norm 1). Thus, by linearity of expectation,
$$ \mathbb{E}_X \sum_{k=1}^{n} \lambda_k \left\langle v_k, X\right\rangle^2 =  \sum_{k=1}^{n} \lambda_k \mathbb{E}_X \left\langle v_k, X\right\rangle^2 = \frac{1}{n} \sum_{k=1}^{n} \lambda_k = \frac{1}{n} (\Delta u)(0).$$
Altogether, we see that
$$  \left( \mathbb{E}_{\omega} \| \omega(t)\|^2 \right)  \cdot   \left( \mathbb{E}_{X} \left\langle X, (D^2 u(0)) X\right\rangle \right)  = t \cdot (\Delta u)(0)$$
and therefore
\begin{align*}
 \mathbb{E} ~u(\omega(t)) &= u(0) +(\Delta u(0))  \frac{t}{2} + \mbox{l.o.t.}
 \end{align*}
Using the Markov property of Brownian motion and letting time run until the Brownian motion impacts the boundary then suggests the representation formula
$$ u(x) = \frac{1}{2} \int_0^{\infty} (-\Delta u)(\omega_x(s)) ds.$$
Since $-\Delta u = f(u)$, we can rewrite this as
$$ u(x) = \mathbb{E} ~\frac{1}{2}\int_0^{\infty}  f(u(\omega_x(s))) ds.$$
A formal justification is, of course, given by the Feynman-Kac formula. This also explains why the expected lifetime of Brownian motion in a domain is given
by the solution of the equation $-\Delta u = 2$: the 2 cancels in the formula.

\begin{proof}[Proof of the Representation Formula] We will now give a proof of the formula for the second derivatives. Our first proof will be stochastic, we will give a sketch of a slightly different argument that is more along the lines of \S 3.1 below.
Using
$$ u(x) = \mathbb{E} ~ \frac{1}{2}\int_0^{\infty}   f(u(\omega_x(s))) ds$$
and the differential quotient
$$ \frac{\partial^2 u}{\partial n^2}(x) = \lim_{\varepsilon \rightarrow 0} \frac{u(x+\varepsilon n) - 2 u(x) + u(x- \varepsilon n)}{\varepsilon^2}.$$
This leads to
$$\frac{\partial^2 u}{\partial n^2}(x) = \lim_{\varepsilon \rightarrow 0} \frac{1}{\varepsilon^2}  \frac{1}{2}~\mathbb{E} \int_0^{\infty} f(u(\omega_{x+\varepsilon n}(s))) - 2 f(u(\omega_{x}(s))) +  f(u(\omega_{x-\varepsilon n}(s))) ds.$$
We note that this is one expectation over three different random processes: a Brownian motion started in $x$ and two more started in $x + \varepsilon n$ and $x - \varepsilon n$. One would assume that these should be related somehow and this can be made precise. We adapt an argument given in \cite{stein2} which is related in spirit. 
 Let the set $A$ denote all Brownian motion paths started in the origin and running for $t$ units of time. Then, for each $x \in \Omega$, we can use translation invariance of Brownian motion in $\mathbb{R}^n$ to write the expectations as an expectation over the set $A$ via
$$  \mathbb{E} \left( f(u(\omega_x(t))) \right) =  \mathbb{E}_{a \in A} \left( f(u(a(t) + x)) \cdot 1_{\left\{a(s) + x \in \Omega ~\mbox{\tiny for all}~0 \leq s \leq t\right\}}\right).$$
This has the advantage of being able to take the expectation with respect to one universal set $A$ shared by all three Brownian motions which differ only by translation invariance. We can use this coupling until the stopping time given by the first time one of the three particles hits the boundary (because it is then absorbed by the boundary and frozen)
$$ \tau = \inf \left\{t >0: \left\{a(t) +x, a(t) +x + \varepsilon n,  a(t) - \varepsilon n \right\} \cap \Omega^c \neq \emptyset \right\}.$$
Until that time, the integral 
$$ X_{\varepsilon} =  \frac{1}{\varepsilon^2} \frac{1}{2} \mathbb{E} _{} \int_0^{\tau} f(u(\omega_{x+\varepsilon n}(s))) - 2 f(u(\omega_{x}(s))) +  f(u(\omega_{x-\varepsilon n}(s))) ds$$
 is actually not difficult to analyze: as $\varepsilon \rightarrow 0$, we can use  $f(u(x)) \in C^2(\Omega)$ and the chain rule
 $$  \frac{\partial^2}{\partial n^2} f(u(x)) = f''(u(x)) \left(\frac{\partial u}{\partial n}\right)^2 + f'(u(x)) \frac{\partial^2 u}{\partial n^2}$$
 then leads to 
 $$ \lim_{\varepsilon \rightarrow 0} X_{\varepsilon} =  \mathbb{E} ~\frac{1}{2} \int_0^{ \tau}  f''(\omega_x(s)) \left( \frac{\partial u}{\partial n}(\omega_x(s)) \right)^2 +  f'(\omega_x(s)) \frac{\partial^2 u}{\partial n^2}(\omega_x(s)) ds.$$

\begin{center}
\begin{figure}[h!]
\begin{tikzpicture}
\draw [very thick] (3,0) to[out=20, in = 270]  (5,1) to[out=90, in =340] (3,2);
\filldraw (5,0.3) circle (0.05cm); 
\node at (6, 1) {$\omega_x(\tau)$};
\node at (4.4, 2) {$\Omega$};
\node at (3.5, 1) {$\Omega^c$};
\node at (2.6, 2.1) {$\partial \Omega$};
\node at (6.4, 1.7) {$\omega_x(\tau)+ \varepsilon n$};
\node at (6.4, 0.3) {$\omega_x(\tau) - \varepsilon n$};
\filldraw (5,1) circle (0.05cm);
\filldraw (5,1.7) circle (0.05cm);
\draw (5, 0.3) -- (5, 1.7);
\end{tikzpicture}
\caption{An example: the middle point has hit the boundary.}
\label{fig:conf}
\end{figure}
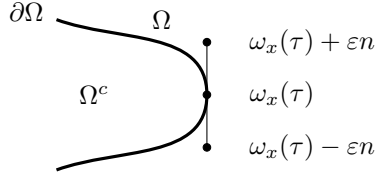
\end{center}

 It remains to control the total contribution of the integral after the first of the three particles has hit the boundary. For the remainder of the argument it will not matter which of the three points has impacted first, the argument is always the same.
 We illustrate the argument using the example shown in Fig. \ref{fig:conf}. In that example, the middle particle has impacted the boundary. 
 Our goal is to understand the remaining integral
 $$ Y_{\varepsilon} =  \frac{1}{\varepsilon^2} \frac{1}{2}\mathbb{E} _{} \int_{\tau}^{\infty} f(u(\omega_{x+\varepsilon n}(s))) - 2 f(u(\omega_{x}(s))) +  f(u(\omega_{x-\varepsilon n}(s))) ds.$$The middle term in the integral has impacted on the boundary and will stay there for all time and not contribute. Thus
  $$ \mathbb{E} _{} \int_{\tau}^{\infty} - 2 f(u(\omega_{x}(s)))   ds = 0 = u(\omega_x(\tau)).$$
As for the two remaining terms, we realize that they correspond to the stochastic representation of the solution and thus
$$  \mathbb{E} _{} ~ \frac{1}{2}\int_{\tau}^{\infty} f(u(\omega_{x+\varepsilon n}(s))) = u(\omega_{x+\varepsilon n}(\tau))$$
$$  \mathbb{E} _{} ~\frac{1}{2}\int_{\tau}^{\infty} f(u(\omega_{x-\varepsilon n}(s))) = u(\omega_{x-\varepsilon n}(\tau)).$$
Altogether, we see that $Y_{\varepsilon}$ is simply the second differential quotient in direction $n$ evaluated at $\omega_x(\tau)$. The harmonic measure
governs the distribution of impact points and from this the desired result follows.
\end{proof}

\begin{proof}[Sketch of Second Proof.] We quickly sketch another proof that is closer to the argument in \S 3.1 and uses the fact that solutions of $-\Delta u = f(u)$ are invariant under translations of the domain $\Omega$. We use the same approach as in \S 3.1 and consider a total of three different domains: $\Omega$ together with $\Omega_+$ and $\Omega_{-}$ (see Fig. 7). We solve $-\Delta u = f(u)$ in each domain and obtain three solutions: $u$, $u_{+}$ and $u_{-}$.
\begin{center}
\begin{figure}[h!]
\begin{tikzpicture}
 \draw[] (0,0) arc (0:360:2cm and 1cm);
   \draw[ultra thick] (0,0) arc (359:348:2cm and 1cm);
  \draw[ultra thick] (0,0) arc (0:140:2cm and 1cm);
  \draw[] (0.2,0.2) arc (0:360:2cm and 1cm);
 \draw[] (0.4,0.4) arc (0:360:2cm and 1cm);
  \draw[ultra thick] (0.4-3.95,0.65) arc (166:322:2cm and 1cm);
  \draw[thick, ->] (1,-0.5) -- (2,0.5);
  \node at (2.2, 0.5) {$n$};
  \node at (-1.6, 0.2) {$ \Omega_{+} \cap \Omega \cap \Omega_{-}$};
\end{tikzpicture}
\caption{Shifting the domain by $\pm \varepsilon$ in direction $n$.}
\end{figure}
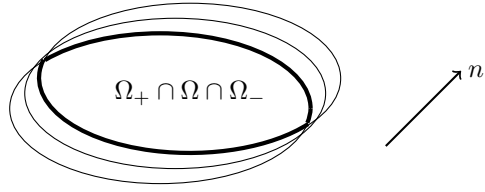
\end{center}
We use again that
$$ \frac{\partial^2 u}{\partial n^2}(x) = \lim_{\varepsilon \rightarrow 0} \frac{u(x+\varepsilon n) - 2 u(x) + u(x- \varepsilon n)}{\varepsilon^2}.$$
Let us fix some small $\varepsilon > 0$ and let us consider the behavior of
$$ w =  \frac{u_{+}(x) - 2 u(x) + u_{-}(x)}{\varepsilon^2} \qquad \mbox{in} ~ \Omega_{+} \cap \Omega \cap \Omega_{-}.$$
We see that, as $\varepsilon \rightarrow 0$, for any $x \in \Omega$,
$$ (-\Delta w)(x) = (1+o(1))  \left[ f''(u(x)) \left(\frac{\partial u}{\partial n}\right)^2 + f'(u(x)) \frac{\partial^2 u}{\partial n^2} \right].$$
Moreover, we also understand the behavior at the boundary since, again as $\varepsilon \rightarrow 0$,
$$ w(x) = (1+o(1)) \frac{\partial^2 u}{\partial n^2} \qquad \mbox{on}~\partial \Omega.$$
We write $w$ as the superposition of two solutions $w = w_1 + w_2$, where $w_1$ has
the correct boundary data while being harmonic
\begin{align*}
\Delta w_1 &= 0 \qquad \quad \mbox{in}~\Omega \\
w_1 &=  \frac{\partial^2 u}{\partial n^2} \qquad \mbox{on}~\partial \Omega
\end{align*}
and where $w_2$ has the correct Laplacian in the domain while satisfying the Dirichlet boundary conditions
\begin{align*}
-\Delta w_2 &=   f''(u(x)) \left(\frac{\partial u}{\partial n}\right)^2 + f'(u(x)) \frac{\partial^2 u}{\partial n^2}\qquad \quad \mbox{in}~\Omega \\
w_2 &=  0 \qquad\qquad\qquad\qquad \qquad \qquad \qquad  \qquad\mbox{on}~\partial \Omega.
\end{align*}
We can represent $w_1$ using harmonic measure
\begin{align*}
w_1(x) &= \int_{\partial \Omega} \frac{\partial^2 u}{\partial n^2}(y) d\omega_x(y)
\end{align*}
and we already discussed the stochastic representation formula for Poisson equations with Dirichlet boundary conditions and thus
\begin{align*}
w_2(x)&= \mathbb{E}~    \frac{1}{2} \int_0^{\infty}  f''(u(\omega_x(s))) \left( \frac{\partial u}{\partial n}(\omega_x(s)) \right)^2 +  f'(u(\omega_x(s))) \frac{\partial^2 u}{\partial n^2}(\omega_x(s)) ds.
\end{align*}
Then $w = w_1 + w_2$ gives the desired representation.
\end{proof}

\subsection{Proof of Theorem 1} 
\begin{proof}Our goal is to show that for each $x \in \Omega$ and each unit norm vector $n$
$$ \frac{\partial^2 u}{\partial n^2}(x) < 0.$$
We recall the representation formula proved in \S 4.2 which reads
\begin{align*}
\frac{\partial^2 u}{\partial n^2}(x) &= \mathbb{E}~  \frac{1}{2} \int_0^{\infty}  f''(u(\omega_x(s))) \left( \frac{\partial u}{\partial n}(\omega_x(s)) \right)^2 + f'(u(\omega_x(s))) \frac{\partial^2 u}{\partial n^2}(\omega_x(s)) ds \\
&+ \int_{\partial \Omega} \frac{\partial^2 u}{\partial n^2}(y) d\omega_x(y).
\end{align*}
We want to argue that the maximum value on the left-hand side,
$$ M = \max_{x \in \Omega \atop \|n\| = 1} \frac{\partial^2 u}{\partial n^2}(x),$$
is negative: $M < 0$. The proof is by contradiction. For this, we use the representation formula in the point $x$ and with the vector $n$ with which $M$ is assumed.
We recall that, by assumption, $f' \geq 0$ and $f'' \leq 0$.
Thus
\begin{align*}
M &\leq  \mathbb{E} \int_0^{\infty}  \frac{1}{2}\cdot f'(u(\omega_x(s))) M ds  + \int_{\partial \Omega} \frac{\partial^2 u}{\partial n^2}(y) d\omega_x(y).
\end{align*}

Since $f'' \leq 0$, we have $0 \leq f'(t) \leq f'(0)$ and
thus
\begin{align*}
 \mathbb{E} \int_0^{\infty}  \frac{1}{2}\cdot f'(u(\omega_x(s))) M ds &\leq M f'(0)  \cdot  \frac{1}{2} \mathbb{E}\int_{0}^{\infty} 1(\omega_x(s)) ds.
 \end{align*}
 Altogether
 $$ M \leq M f'(0)  \cdot  \frac{1}{2} \mathbb{E}\int_{0}^{\infty} 1(\omega_x(s)) ds +  \underbrace{\int_{\partial \Omega} \frac{\partial^2 u}{\partial n^2}(y) d\omega_x(y)}_{< 0}.$$
In order for this inequality to be satisfied for some $M > 0$, we require
$$  f'(0)  \cdot  \frac{1}{2} \mathbb{E} \int_{0}^{\infty} 1(\omega_x(s)) ds > 1.$$
We will now prove an upper bound on the integral that will lead to a contradiction.
It is know that among all domains with fixed measure, the maximum lifetime is maximized for a
ball and
can be bounded from above by the maximum of the solution of $-\Delta \psi = 2$ in a ball with the same radius as $\Omega$. The solution is given by
$$ \psi(x_1, \dots, x_n) = \frac{r^2}{n} - \frac{x_1^2 -x_2^2 - \dots - x_n^2}{n}$$
and thus
$$ \max_{x \in \Omega} ~\mathbb{E} \int_{0}^{\infty} 1(\omega_x(s)) ds \leq \frac{r^2}{n} = \frac{ |\Omega|^{\frac{2}{n}}}{n \omega_n^{2/n}}.$$
As noted, it would suffice to show that
$$ \max_{x \in \Omega}~ f'(0)  \cdot  \frac{1}{2} \mathbb{E}\int_{0}^{\infty} 1(\omega_x(s)) ds \leq 1,$$
which, using the estimate on the expected lifetime, follows as soon as
 $$  f'(0) \leq \frac{2 n \omega_n^{2/n}}{|\Omega|^{2/n}}.$$
\end{proof}

\subsection{Proof of Theorem 2}
\begin{proof}
Let 
$$ -\Delta u = f(u) \qquad \mbox{in}~\Omega.$$
We can then consider the solution (subject to Dirichlet boundary conditions) of
$$ -\Delta v = f(u)^* \quad \mbox{in}~B,$$
where $B$ is the ball centered at the origin having the same measure as $\Omega$, $|B| = |\Omega|$, and $f(u)^*$ denotes the symmetric decreasing rearrangement of $f(u)$ (we refer to the book of Baernstein \cite{baern} for a detailed introduction to symmetrization techniques). By Talenti's rearrangement principle \cite{talenti, talenti2}, we have
$$ v \geq u^*$$ 
and, in particular, we have $\| v\|_{L^{\infty}(B)} \geq \| u\|_{L^{\infty}(\Omega)}$. We will compare this to the solution (subject to Dirichlet boundary conditions) of 
$$ -\Delta \psi = f(\psi) \quad \mbox{in}~B.$$
We will show that $v(x) \leq \psi(x)$ for all $x \in B$.
Observe that both $v$ and $\psi$ are nonnegative functions. Moreover, since $v \geq u^*$ and $f$ is non-decreasing, 
$$ -\Delta v =f(u)^* =  f(u^*) \leq f(v).$$
We define the operator $(-\Delta)^{-1} f$ as mapping a given function $f:B \rightarrow \mathbb{R}_{\geq 0}$ to the solution $w:B \rightarrow \mathbb{R}$ of $-\Delta w = f$ with Dirichlet boundary conditions. The maximum principle implies that if, for all $x \in B$, we have $g_1(x) \leq g_2(x)$, then, for all $x \in B$,
$$ \left[(-\Delta)^{-1} g_1\right](x) \leq  \left[(-\Delta)^{-1} g_2\right](x).$$

 From this we obtain, for all $x \in \Omega$, the inequality $ v(x) \leq \left[(-\Delta)^{-1} f(v)\right](x)$.
 Our goal is to show that, for all $x \in B$
 $$ \psi(x) \geq v(x).$$
Assume the statement is false and that, for $x_0 \in B$, we have $v(x_0) > \psi(x_0)$. Then
\begin{align*}
0 < v(x_0) - \psi(x_0) &\leq \left[ (-\Delta)^{-1} f(v) - \psi \right](x_0) \\
 &= \left[ (-\Delta)^{-1} \left( f(v) - f(\psi) \right) \right](x_0) \\
 &\leq \left[ (-\Delta)^{-1}  \max_{y \in B}  \left[ f(v(y)) - f(\psi(y))\right]  \right](x_0) \\
  &\leq \max_{y \in B} \left[ f(v(y)) - f(\psi(y)) \right] \cdot \left[ (-\Delta)^{-1}~ 1  \right](x_0).
 \end{align*}
 The mean-value theorem yields
 $$  \max_{y \in B} \left[ f(v(y)) - f(\psi(y)) \right]  \leq \left( \max_{t> 0} f'(t) \right) \cdot  \max_{y \in B} \left[ v(y) - \psi(y) \right].$$
 The function $(-\Delta)^{-1}~1$ is merely the solution of $-\Delta u = 1$ and well understood. In the proof of Theorem 3, we showed the inequality
 $$ \left[ (-\Delta)^{-1}~ 1  \right](x) \leq \frac{ |\Omega|^{\frac{2}{n}}}{2n \omega_n^{2/n}}.$$
Hence,
$$  \left(v - \psi\right)(x_0)\leq    \frac{ |\Omega|^{\frac{2}{n}}}{2n \omega_n^{2/n}} f'(0) \cdot \max_{y \in B} \left[ v(y) - \psi(y) \right].$$
However, since $x_0$ was arbitrary (as long as $v(x_0) > \psi(x_0)$), we may as well pick $x_0$ to be the point in which $v(x) - \psi(x)$
assumes its maximum. This leads to a contradiction as soon as
 $$\frac{ |\Omega|^{\frac{2}{n}}}{2n \omega_n^{2/n}} \cdot f'(0) < 1.$$
 This contradiction then shows that
 $$ \|u\|_{L^{\infty}(\Omega)} =  \|u^*\|_{L^{\infty}(B)} \leq  \|v\|_{L^{\infty}(\Omega)}  \leq  \|\psi\|_{L^{\infty}(\Omega)}.$$
\end{proof}


\begin{thebibliography}{}

\bibitem{arango} J. Arango, A. Gomez,
Critical points of solutions to elliptic problems in planar domains. 
Commun. Pure Appl. Anal. 10 (2011), no. 1, 327--338.

\bibitem{baern} A. Baernstein, 
Symmetrization in analysis.
With David Drasin and Richard S. Laugesen. With a foreword by Walter Hayman. New Mathematical Monographs, 36. Cambridge University Press, Cambridge, 2019.

\bibitem{bandle} C. Bandle, 
On isoperimetric gradient bounds for Poisson problems and problems of torsional creep.
Z. Angew. Math. Phys. 30 (1979), no. 4, 713--715.


\bibitem{banuelos} R. Banuelos, R. Latala and P. Mendez-Hernandez, A Brascamp-Lieb-Luttinger-type
inequality and applications to symmetric stable processes, Proc. Amer. Math. Soc. 129
(2001), 2997--3008.


\bibitem{beck1} T. Beck, Uniform level set estimates for ground state eigenfunctions. SIAM J. Math. Anal. 50 (2018), no. 4, 4483--4502.

\bibitem{beck2} T. Beck, The torsion function of convex domains of high eccentricity, Potential Analysis 53, p. 701--726 (2020).

\bibitem{beck3} T. Beck, Localization of the first eigenfunction of a convex domain, Comm. PDE, to appear.

\bibitem{bian} B. Bian and P. Guan, A microscopic convexity principle for nonlinear partial differential
equations, Invent. Math. 177 (2009), p. 307--335

\bibitem{borell} C. Borell, Brownian motion in a convex ring and quasiconcavity, Comm. Math. Phys. 86 (1982), p. 143--147.

\bibitem{borell2} C. Borell, A note on parabolic convexity and heat conduction, Annales de l'I.H.P. Probabilités et statistiques, Tome 32 (1996) p. 387--393.

\bibitem{brascamp} H. J. Brascamp and E. H. Lieb, On extensions of the Brunn-Minkowski and Prekopa Leindler theorems, including inequalities for log concave functions, and with an application to the diffusion
equation, J. Funct. Anal. 22 (1976), 366--389.

\bibitem{caf} L. Caffarelli and J. Spruck, Convexity properties of solutions to some classical variational problems, Comm. Partial Differ. Equations, 7 (1982), 1337--1379.


\bibitem{caf2} L. Caffarelli L. and A. Friedman, Convexity of solutions of some semilinear elliptic equations,
Duke Math. J., 52(1985), 431--455.


\bibitem{cha} S.-Y. A. Chang, X.-N. Ma and P. Yang, Principal curvature estimates for the convex level sets of semilinear elliptic equations. Discrete Contin. Dyn. Syst. 28 (2010), no. 3, 1151--1164. 

\bibitem{chau1} A. Chau and B. Weinkove, The Stefan problem and concavity, arXiv:2004.04284 

\bibitem{finn} D. Finn,
Convexity of level curves for solutions to semilinear elliptic equations. 
Commun. Pure Appl. Anal. 7 (2008), no. 6, 1335--1343.


\bibitem{gab} R. Gabriel, A result concerning convex level surfaces of 3-dimensional harmonic functions, J. London Math.Soc., 32 (1957), 286--294.


\bibitem{gidas}  B. Gidas, W. Ni, and L. Nirenberg. Symmetry and related properties via the maximum
principle. Comm. Math. Phys., 68: p. 209--243, 1979.

\bibitem{gugg} H. Guggenheimer, 
Concave solutions of a Dirichlet problem.
Proc. Amer. Math. Soc. 40 (1973), p. 501--506

\bibitem{hamel} F. Hamel, N. Nadirashvili and Y. Sire, Yannick
Convexity of level sets for elliptic problems in convex domains or convex rings: two counterexamples. 
Amer. J. Math. 138 (2016), no. 2, 499--527. 

\bibitem{hen} Antoine Henrot, Ilaria Lucardesi and Gerard Philippin, On two functionals involving the maximum of the torsion function, ESAIM: COCV 24 (2018), p. 1585--1604


\bibitem{henrot} A. Henrot, C. Nitsch, P. Salani \& C. Trombetti,
Optimal Concavity of the Torsion Function,
Journal of Optimization Theory and Applications 178 (2018), p. 26--35

\bibitem{hoskins} J. Hoskins and S. Steinerberger, Towards Optimal Gradient Bounds for the Torsion Function in the Plane, arXiv:1912.08376

\bibitem{kawohl0} B. Kawohl, Rearrangements and convexity of level sets in PDE. Lecture Notes in Mathematics, 1150. Springer-Verlag, Berlin, 1985.

\bibitem{kawohl1} B. Kawohl, Variations on the p-Laplacian. Nonlinear elliptic partial differential equations, 35--46, Contemp. Math., 540, Amer. Math. Soc., Providence, RI, 2011.

\bibitem{keady0} G. Keady, The power concavity of solutions of some semilinear elliptic boundary-value problems, Bull. Aust. Math. Soc. 31 (1985), 181--184

\bibitem{keady} G. Keady and A. McNabb, 
The elastic torsion problem: solutions in convex domains. 
New Zealand J. Math. 22 (1993), no. 2, 43--64.

\bibitem{ken} A. U. Kennington, Power concavity and boundary value problems, Indiana Univ. Math. J. 34 3 (1985), 687--704

\bibitem{kor} N. Korevaar, Convex solutions to nonlinear elliptic and parabolic boundary value problems. Indiana Univ. Math J. 32 (1983), p. 603--614.

\bibitem{kor2} N. Korevaar and J. L. Lewis,
Convex solutions of certain elliptic equations have constant rank hessians,
Archive for Rational Mechanics and Analysis  97, p. 19--32 (1987)

\bibitem{kos0} A. A. Kosmodem’yanskii, Sufficient conditions for the concavity of the solution
of the Dirichlet problem for the equation $\Delta u =- 1$, Math. Notes of Acad. Sci.
of U.S.S.R. 42 (1987) 798-801.

\bibitem{stv1} A. Kosmodemyanskii,
The behavior of solutions of the equation $\Delta u = -1$ in convex domains. (Russian)
Dokl. Akad. Nauk SSSR 304 (1989), no. 3, 546--548; translation in
Soviet Math. Dokl. 39 (1989), 112--114


\bibitem{lions} P.-L. Lions, Two geometrical properties of solutions of semilinear problems. 
Applicable Anal. 12 (1981), 267--272. 


\bibitem{lu} J. Lu and S. Steinerberger, Optimal Trapping for Brownian Motion: a Nonlinear Analogue of the Torsion Function
Potential Analysis 54, p. 687--698 (2021).

\bibitem{ma} X-N, Ma, S. Shi and Y. Ye, 
The convexity estimates for the solutions of two elliptic equations. 
Comm. Partial Differential Equations 37 (2012), no. 12, 2116--2137.

\bibitem{makar} L. G. Makar-Limanov, 
The solution of the Dirichlet problem for the equation $\Delta u = -1$ in a convex region. 
Mat. Zametki 9 (1971), 89--92.

\bibitem{maxpayne} L. E. Payne, Bounds for the maximum stress in the Saint Venant torsion problem, Indian J. Mech. Math (1968), p. 51--59.

\bibitem{pol} G. Polya, Torsional rigidity, principal frequency, elecostactic capacity and symmetrization, Quart. Appl. Math. 6, p. 267 -- 277 (1948).

\bibitem{reilly} R. Reilly, Mean Curvature, The Laplacian, and Soap Bubbles,  The American Mathematical Monthly 89 (1982),
p. 180--198

\bibitem{paynestak} I. Stakgold and L. E. Payne. Nonlinear problems in nuclear reactor analysis. Nonlinear problems in the physical sciences and biology. Springer, Berlin, Heidelberg, 1973. 298--307.



\bibitem{sperb} R. Sperb, Maximum principles and their applications,  Mathematics in
Science and Engineering, vol. 157, Academic Press, New York, 1981.


\bibitem{stein} S. Steinerberger, Topological Bounds on Fourier Coefficients and Applications to Torsion, Journal of Functional Analysis 274 (2018),  p. 1611-1630.

\bibitem{stein2} S. Steinerberger, A Pointwise Inequality for Derivatives of Solutions of the Heat Equation in Bounded Domains, arXiv:2102.02736

\bibitem{talenti}  G. Talenti, Elliptic equations and rearrangements. Ann. Scuola Norm. Sup. Pisa Cl. Sci. (4) 3 (1976), no. 4, 697--718. 

\bibitem{talenti2} G. Talenti. Nonlinear elliptic equations, rearrangements of functions and Orlicz spaces.
Ann. Mat. Pura Appl., 120 (1977): p. 159--184.

\bibitem{weitsman} A. Weitsman, Symmetrization and the Poincare metric. Ann. of Math. 124 (1986), p.159–169.

\end{thebibliography}
\end{document}